\documentclass{article}

\usepackage{amsmath}
\usepackage{amssymb}
\usepackage{latexsym}
\usepackage{amsxtra}
\usepackage{amscd}
\usepackage{theorem}
\usepackage{xypic}

\theoremstyle{plain}

{\theorembodyfont{\slshape}

        \newtheorem{thm}{Theorem}[section]
        \newtheorem{cor}[thm]{Corollary}
        \newtheorem{lem}[thm]{Lemma}
        \newtheorem{prop}[thm]{Proposition}

}

{\theorembodyfont{\rmfamily}

        \newtheorem{defn}[thm]{Definition}
        \newtheorem{prob}[thm]{Problem}
        \newtheorem{rem}[thm]{Remark}
        
        \newtheorem{exa}[thm]{Example}

        \newtheorem{question}[thm]{Question}

}

\renewcommand{\em}{\sl}

\newcommand{\proof}{{\bf Proof:\ }}
\newcommand{\Endproof}{\hspace*{\fill} $\Box$ \vspace{1ex} \noindent }

\makeatletter
\renewcommand{\subsection}{\@startsection{subsection}{2}%
        {\z@}{-3.25ex plus -1ex minus-.2ex}{-1em}{\bf}}
\makeatother

\newcommand{\ZZ}{\mathbb{Z}}
\newcommand{\QQ}{\mathbb{Q}}

\newcommand{\NN}{\mathbb{N}}

\newcommand{\OO}{\mathcal{O}}
\newcommand{\I}{\mathcal{I}}

\newcommand{\Aut}{{\rm Aut}}

\newcommand{\Gal}{{\rm Gal}}
\newcommand{\Ker}{{\rm Ker}}

\newcommand{\Tr}{{\rm Tr}}

\newcommand{\GL}{{\rm GL}}

\newcommand{\m}{\mathfrak{m}}
\newcommand{\p}{\mathfrak{p}}
\newcommand{\Spec}{{\rm Spec\,}}

\newcommand{\Frac}{{\rm Frac}}

\newcommand{\rank}{{\rm rank}}
\newcommand{\Der}{{\rm Der}}

\newcommand{\inj}{\hookrightarrow}
\newcommand{\To}{\;\longrightarrow\;}
\newcommand{\iso}{\stackrel{\sim}{\to}}

\newcommand{\lpfeil}[1]{\stackrel{#1}{\To}}
\newcommand{\gen}[1]{\mathopen\langle#1\mathclose\rangle}

\newcommand{\Kb}{\bar{K}}
\newcommand{\Lb}{\bar{L}}
\newcommand{\Trb}{\overline{\Tr}}
\newcommand{\Bb}{\bar{B}}
\newcommand{\Ab}{\bar{A}}
\newcommand{\xb}{\bar{x}}
\newcommand{\Bh}{\hat{B}}
\newcommand{\Ah}{\hat{A}}

\begin{document}

\title{Regularity of quotients by an automorphism of order $p$} 
\author{Stefan Wewers\\ IAZD, Leibniz-Universit\"at Hannover}
\date{}

\maketitle

\begin{abstract}
  Let $B$ be a regular local ring and $G\subset\Aut(B)$ a finite group of
  local automorphisms. Assume that $G$ is cyclic of prime order $p$, where $p$
  is equal to the residue characteristic of $B$. We give conditions under
  which the ring of invariants $A=B^G$ is again regular. 
\end{abstract}

\section*{Introduction}

Let $X$ be a regular integral scheme and $G\subset\Aut(X)$ a finite group of
automorphisms of $X$. The quotient scheme $Y:=X/G$ may not be regular;
its singularities are, by definition, {\em quotient singularities}. To study
the singularities of $Y$, we may localize and assume that $X=\Spec B$ and
$Y=\Spec A$, where $B$ is a local domain and $A=B^G$ is the ring of
invariants.  

Quotient singularities have been intensively studied, in connection with
reso\-lution of singularities and as objects in their own right. However, most
results concern {\em tame} quotient singularities. In the above situation this
means that the order of $G$ is prime to the characteristic of the residue
field of $B$. 

In a recent preprint \cite{Lorenzini09}, D.\ Lorenzini has studied the
resolution graphs of {\em wild} quotient singularities in dimension $2$,
exhibiting many interesting features that do not occur for tame quotient
singularities. His results rely heavily on a detailed combinatorial study of
the possible intersection matrices that can occur. The present note aims at
complementing the methods used in \cite{Lorenzini09}. The basic idea is the
following.

Let $Y=X/G$ be as above, and let $Y'\to Y$ be a resolution of $Y$. Then the
fiber product $X':=X\times_Y Y'$ is a $G$-equivariant modification of $X$
whose quotient $Y'=X'/G$ is regular. Conversely, given a $G$-equivariant
modification $X'\to X$, the quotient scheme $Y':=X'/G$ is a modification of
$Y$ -- which may or may not be regular. From this point of view it is natural
to look for conditions on $X'$ under which $Y'$ is regular. However, it is
difficult to find such conditions in the literature which apply to the
case of wild group actions.




Our main result (Theorem \ref{mainthm}) gives a sufficient condition for the
regularity of the quotient scheme $Y=X/G$ when $G$ is a cyclic group of order
$p$. This is admittedly a very modest contribution to our motivating
problem. Still, our criterion seems to be new, and we hope that it will be
useful in the future.

The motivation for writing this paper grew out of discussions with F.\
Kir\'aly, who discovered a special case of our main result, see Example
\ref{exa1.2}. The author thanks him, W.\ L\"utkebohmert and M.\ Raynaud for
helpful discussions and comments on earlier versions of this paper.

\section{The main result} \label{sec1}

\subsection{} \label{sec1.1}

Let $(B,\m,k)$ be a noetherian regular local domain, and let $G\subset\Aut(B)$
be a finite group of local automorphisms. We are interested in
the following question.

\begin{prob} \label{mainprob}
 Under which condition is the ring of invariants
 \[
       A :=B^G = \{\, x\in B \mid 
          \sigma(x)=x\;\;\text{for all $\sigma\in G$} \,\}
 \]
 again regular?
\end{prob}

For the applications we have in mind, the following additional assumptions
seem reasonable and useful:
\begin{itemize}
\item[(a)]
  The residue field $k$ is algebraically closed.
\item[(b)]
  The induced $G$-action on $k$ is trivial.
\item[(c)]
  $A$ is noetherian and $B$ is a finite $A$-algebra.
\end{itemize}
Assumptions (a) and (b) can always be achieved by passing to the strict
henselization and are therefore harmless. They imply that $A$ is a local domain
with residue field $k$. Condition (c) is satisfied in most situations
arising from a geo\-metric context. Suppose, for instance, that $B$ is the
localization of a finitely generated algebra over an excellent domain $R$, and
that the action of $G$ on $B$ is $R$-linear. Then (c) holds. See
\cite{Fogarty80} for a general discussion of Condition (c). 

If the order of the group $G$ is prime to the characteristic of the residue
field $k$, a definitive answer to Problem \ref{mainprob} is known.

\begin{thm} \label{serrethm} 
  Suppose that (a)-(c) holds and that the order of
  $G$ is prime to the characteristic of the residue field $k$. Then the ring
  $A$ is regular if and only if the image of $G$ in $\GL(\m/\m^2)$ is
  generated by {\em pseudo-reflections}\footnote{An element $\sigma\in \GL(V)$
    is called a {\em pseudo-reflection} if $\rank(\sigma-1)\leq 1$}.
\end{thm}

\proof
See \cite{Serre67} or \cite{Watanabe76}
\Endproof

The main result of the present paper (Theorem \ref{mainthm} below) gives a
sufficient criterion for $A$ to be regular in the case where $G$ is cyclic of
prime order $p$. Here $p$ may be equal to the characteristic
of $k$, and hence our result is not covered by Theorem \ref{serrethm}

\subsection{} \label{sec1.2}

Let us now assume that $G$ is cyclic of order $p$, where $p$ is
prime. We choose a generator $\sigma\in G$ and consider the ideal
\[
       \I_\sigma :=\gen{\sigma(x)-x \mid x\in B}_B\lhd B.
\]
By definition, $\I_\sigma$ is the smallest $G$-invariant ideal such that $G$
acts trivially on $B/\I_\sigma$. In particular, this shows that $\I_\sigma$
does not depend on the chosen generator $\sigma$. Condition (b) says that
$\I_\sigma$ is contained in $\m$. 

An element $x\in B$ is called a {\em regular parameter} if
$x\in\m\backslash\m^2$. Since $B$ is regular, this means that $x$ is part of a
regular system of parameters. It follows that $B/Bx$ is a regular local ring
and that $(x)\lhd B$ is a prime ideal.

Here is our main result:

\begin{thm} \label{mainthm} Suppose that there exists a regular parameter
  $\pi\in\m\backslash\m^2$ such that 
  \begin{enumerate}
  \item
    $\I_\sigma=(\pi^\delta)$, with $\delta\in\NN$, 
  \item
    we either have $\I_\sigma=(\sigma(\pi)-\pi)$, or else 
    $\sigma(\pi)-\pi\in(\pi^{\delta+1})$.
  \end{enumerate} 
  Then $A=B^G$ is regular.
\end{thm}
 
After some preliminary work done in Section \ref{sec2} and \ref{sec3}, the
proof of Theorem \ref{mainthm} will be given in Section \ref{sec4}. For the
rest of this section, we discuss the scope and the limitations of Theorem
\ref{mainthm} and some open problems.

\begin{rem} \label{rem1.1} Let $\bar{\sigma}\in GL(\m/\m^2)$ denote the image
  of $\sigma$. By definition of $\I_\sigma$, $\m/(\I_\sigma+\m^2)$ is the
  largest quotient of $\m/\m^2$ on which $\bar{\sigma}$ acts trivially. If the
  hypothesis of Theorem \ref{mainthm} is verified, then it follows that
  $\bar{\sigma}$ is a pseudo-reflection. So for ${\rm char}(k)\neq p$, Theorem
  \ref{mainthm} is a direct consequence of Theorem \ref{serrethm}. Moreover,
  in this case the natural homomorphism
  \[
        G \to \GL(\m/\m^2)
  \]
  is known to be injective (see e.g.\ \cite{Watanabe76}). This implies that
  we must have $\delta=1$ in Condition (i) of Theorem \ref{mainthm} and that
  $\I_\sigma=(\sigma(\pi)-\pi)$ in Condition (ii). 

  On the other hand, if ${\rm char}(k)=p$, $\delta$ may be strictly larger than
  $1$. In particular, the hypothesis of Theorem \ref{mainthm} is stronger
  then the condition `$\bar{\sigma}$ is a pseudo-reflection': it is easy to
  give examples where $\bar{\sigma}=1$ and $A$ is not regular. See e.g.\
  \cite{Peskin83}. 
\end{rem}

\begin{rem}
  The author suspects that Condition (ii) in Theorem \ref{mainthm} is implied
  by Condition (i) and could therefore be omitted. However, he did not succeed
  in proving this.
\end{rem} 

\begin{rem}
  The sufficient condition given by Theorem \ref{mainthm} is not necessary for
  $A$ to be regular. See e.g.\ Example \ref{exa1.1} below. 
\end{rem}

The Purity Theorem of Zariski-Nagata gives us a necessary condition for
regularity. Namely, if $A$ is regular, then all minimal prime ideals of $B$
containing $\I_\sigma$ have height one, and are therefore principal. 
This prompts the following question:

\begin{question} \label{question1.1}
  Is it true that $A$ is regular if and only if $\I_\sigma$ is a principal
  ideal?
\end{question}

\subsection{}

We shall give three examples that illustrate various points. In all three
examples, the ring $B$ is a ring of power series over a complete discrete
valuation ring $R$. In particular, $B$ has dimension $2$. We assume moreover
that the action of the group $G$ fixes the subring $R\subset B$.  

The first example shows that the hypothesis of Theorem \ref{mainthm} is
not a necessary condition for regularity of the ring $A$. 

\begin{exa} \label{exa1.1}
  Let $R$ be a complete discrete valuation ring of mixed characteristic
  $(0,p)$, containing a $p$th root of unity $\zeta_p$. Set $B:=R[[x]]$, and
  let $\sigma:B\iso B$ be the $R$-automorphism of order $p$ given by
  \[
        \sigma(x) = \zeta_p\cdot x.
  \]
  Then 
  \[
      A:=B^{\gen{\sigma}} = R[[x^p]]
  \]
  is regular. However, the ideal
  \[
      \I_\sigma=\big(\,(\zeta_p-1)x\,\big)\lhd B
  \]
  is contained in two distinct principal prime ideals, so Condition (i) of
  Theorem \ref{mainthm} fails.
\end{exa} 

\begin{rem} \label{rem1.4}
  Example \ref{exa1.1} is a special case of the following situation. Assume
  that the ring $B$ has a regular system of parameters
  $(\pi_1,\ldots,\pi_{d-1},x)$ such that
  \[
      \sigma(\pi_i)=\pi_i \quad\text{for $i=1,\ldots,d-1$}
  \]
  and
  \[
     \sigma(x) \equiv u\cdot x \mod{(\pi_1,\ldots,\pi_{d-1})}
  \]
  for a unit $u\in B^\times$. Then $A=B^G$ is regular, see \cite{KatzMazur},
  Proposition 7.5.2. In the situation of Theorem \ref{mainthm} a system
  of parameters $(\pi_1,\ldots,\pi_{d-1},x)$ as above exists, but this is a
  consequence of the proof of Theorem \ref{mainthm}, and is not obvious
  beforehand. 
\end{rem}

The next example, which was first studied by F.\ Kir\'aly
\cite{FranzDiss}, describes a special situation where the answer to 
Question \ref{question1.1} is affirmative. 

\begin{exa} \label{exa1.2}
  Let $R$ and $B=R[[x]]$ be as in Example \ref{exa1.1}. Let $\sigma:B\iso B$
  be an automorphism of order $p$ which induces the identity on the residue
  field $k$. Contrary to Example \ref{exa1.1}, we assume
  that $\sigma$ restricts to a nontrivial automorphism of the subring
  $R\subset B$.

  Let $v:R\to\ZZ\cup\{\infty\}$ denote the discrete valuation on $R$ and let
  $\pi\in R$ be a uniformizer, i.e.\ $v(\pi)=1$. Our assumptions imply that
  \[
        \sigma(\pi)= \pi + u\pi^\mu,
  \]
  with $\mu\geq 2$ and $u\in R^\times$. Write 
  \[
      f_\sigma:=\sigma(x)-x = a_0+a_1\,x+a_2\,x^2+\ldots,
  \]
  with $a_i\in R$ and set
  \[
       \nu:=\min\{\,v(a_i) \mid i=0,1,\ldots\}.
  \]
  Since $(\pi,x)$ is a system of parameters for $B$, we have
  \[
      \I_\sigma=(\pi^\mu,f_\sigma)\lhd B.
  \]
  It follows that $\I_\sigma$ is a principal ideal if and only if one
  of the following cases occurs:
  \begin{itemize}
  \item[(I)]
    $\mu\leq \nu$, or
  \item[(II)]
    $v(a_0)=\nu<\mu$.
  \end{itemize} 
  In Case (I) we have $\I_\sigma=(\pi^\mu)$ and in Case (II) we have
  $\I_\sigma=(\pi^\nu)$. In both cases, the hypothesis of Theorem
  \ref{mainthm} holds, and hence $A=B^G$ is regular. In this special case, the
  statement of Theorem \ref{mainthm} has been proved earlier by F.\ Kir\'aly,
  and his results yield somewhat more. In Case (I), the ring of invariants $A$
  is actually a power series ring over $R^G$. Moreover, if neither Case
  (I) nor Case (II) holds, then $\I_\sigma$ is not a principal ideal and the
  ring $A$ is not regular. We refer to \cite{FranzDiss} for more details. 
\end{exa}

The distinction of Case (I) and (II) in Example \ref{exa1.2} illustrates a
dichotomy which is crucial for the proof of Theorem \ref{mainthm}. Let
$L:=\Frac(B)$ and $K:=\Frac(A)$ denote the fraction fields. Then $L/K$ is a
Galois extension with Galois group $G$. Assume that the hypothesis of Theorem
\ref{mainthm} holds. Let $v$ denote the discrete valuation of $L$ which
corresponds to the localization of $B$ at the prime ideal $(\pi)$. It follows
from Condition (i) that $v$ is ramified in $L/K$. The two cases of Condition
(ii) correspond to Case (I) and Case (II) in Example \ref{exa1.2}. We shall
see in Section \ref{sec3} that in the first case the extension $L/K$ is {\em
  totally ramified} along $v$, whereas it is {\em fiercely ramified} in the
second case (see Definition \ref{def3.1} and Proposition \ref{prop3.1}).

\subsection{}

If we drop the assumption that the group $G$ acting on $B$ is cyclic of order
$p$, the study of the ramification behavior of the extension $L/K$ becomes
much more complicated. The problem is that it is in general difficult to
`separate' a wildly ramified extension in a canonical way into subextensions
which are either totally or fiercely ramified. This seems to be the main
reason why the statement of Theorem \ref{mainthm} does not easily generalize
to more general groups. Our last example illustrates this point.

\begin{exa} \label{exa1.3} Let $R$ and $B=R[[x]]$ be as in Example
  \ref{exa1.1} and \ref{exa1.2}. Let $G\subset\Aut(B)$ be a finite group of
  automorphisms which fixes the subring $R\subset B$ and acts trivially on the
  residue field $k$. We assume that $G$ is an elementary abelian group of
  order $p^2$, i.e.\ $G\cong\ZZ/p\times\ZZ/p$, with generators
  $\sigma_1,\sigma_2$. We also assume that the induced action of $G$ on $R$ is
  faithful. Let $A_i:=B^{\gen{\sigma_i}}$, $i=1,2$. The restriction of
  $\sigma_1$ to $A_2$ is an automorphism of order $p$ such
  that $A_2^{\sigma_1}=A$. Likewise, $A_1^{\sigma_2}=A$.

  Suppose that $A_1$ and $A_2$ are regular. Does it follow that $A$ is
  regular?  The statement of Theorem \ref{serrethm} may lead one to believe
  that the answer to this question is yes. However, this is not the case; in
  the following we shall give an explicit counterexample, where $p=2$.

  Let $K:=\QQ_2^{\rm nr}$ denote the maximal unramified extension of
  $\QQ_2$. Let $L/K$ be the Galois extension generated by the roots of the
  polynomial $(x^2-2)(x^2-3)$. Hence $L/K$ is generated by elements
  $\sqrt{2},\sqrt{3}$ which are square roots of $2$ and $3$, respectively. Let
  $G=\Gal(L/K)$ denote the Galois group. Then $G\cong\ZZ/2\times\ZZ/2$ is
  generated by the two automorphisms $\sigma_1,\sigma_2$ determined by
  \[\begin{split}
     \sigma_1(\sqrt{2})=\sqrt{2}, \quad & \sigma_1(\sqrt{3})=-\sqrt{3},\\
     \sigma_2(\sqrt{2})=-\sqrt{2}, \quad &\sigma_2(\sqrt{3})=\sqrt{3}.
  \end{split}\]
  It follows that
  \[
      K_1:=\QQ_2(\sqrt{2})=L^{\sigma_1}, \quad 
      K_2:=\QQ_2(\sqrt{3})=L^{\sigma_2}.
  \] 
  Let $v$ denote the unique extension to $L$ of the $2$-adic valuation,
  normalized by $v(2)=1$. Set
  \[
      \pi_1:=\sqrt{2}, \quad \pi_2:=\sqrt{3}-1,\quad
      \pi:=\frac{\pi_2}{\pi_1}-1.
  \]
  One checks that
  \[
       v(\pi_1)=v(\pi_2)=\frac{1}{2}, \quad v(\pi)=\frac{1}{4}.
  \]
  It follows that
  \[
     \OO_{K_1}=\ZZ_2^{\rm nr}[\pi_1], \quad 
     \OO_{K_2}=\ZZ_2^{\rm nr}[\pi_2], \quad
     \OO_K=\ZZ_2^{\rm nr}[\pi].
  \]

  Let $B:=\OO_L[[x]]$ be the ring of power series in $\OO_L$. We extend the
  action of $G$ on $\OO_L$ to $B$ by setting
  \[
        \sigma_1(x):=x,\qquad \sigma_2(x)=-x.
  \] 
  This is the  example announced above. Indeed, it is easy to see that 
  \[
      A_1=B^{\sigma_1}=\OO_{K_1}[[x]]  
  \]
  is a power series ring over $\OO_{K_1}$ and hence regular. Similarly, if we
  set $y:=(1+\pi)x$, then we find that
  \[
        B=\OO_L[[y]], \quad \text{and}\quad \sigma_2(y)=y.
  \]
  It follows that
  \[
      A_2=B^{\sigma_2}=\OO_{K_2}[[y]]  
  \]
  is regular, too. However, $A=B^G=A_1^{\sigma_2}$ is not regular. This can be
  checked using the if-and-only-if criterion from Example \ref{exa1.2}. Note
  that $(\pi_1,x)$ is a regular sequence of parameters for $A_1$. Now
  \[
       \sigma_2(\pi_1)-\pi_1=-2\pi_1, \quad \sigma_2(x)-x= -2x,
  \]
  and so
  \[
     \I_{\sigma_2}=(-2\pi_1,-2x) = (2)\cdot\m_{A_1} \lhd A_1
  \]
  is not a principal ideal. Using the criterion of \cite{FranzDiss} mentioned
  in Example \ref{exa1.2}, we conclude that $A$ is not regular. The reader is
  invited to check this by a direct calculation of the dimension of
  $\m_A/\m_A^2$. 
\end{exa}

\section{Derivations and $p$-cyclic inseparable descent} \label{sec2}

In this section we prove an auxiliary result (Corollary \ref{cor2.1}) which is
a crucial step in the proof of Theorem \ref{mainthm}. The setup is similar as
in the previous section; however, we work in equal characteristic $p$, and
instead of an automorphism of order $p$ we consider a derivation of the ring
$B$. 

Let $(B,\m,k)$ denote a noetherian regular local ring of
dimension $d\geq 1$. We assume moreover that $B$ is complete and has
characteristic $p$. Then it follows from \cite{MatsumuraCA}, Corollary 2 of
Theorem 60, that $B$ is isomorphic to a formal power series ring over
$k$. More precisely, if $(x_1,\ldots,x_d)$ is a regular system of parameters
for $B$, then we get an isomorphism
\[
        B\cong k[[x_1,\ldots,x_d]].
\]
We let $L$ denote the fraction field of $B$. We write $\Der_k(B)$ for the
$p$-Lie-algebra of $k$-derivations $\theta:B\to B$. Note that any
such derivation extends uniquely to a (continuous) derivation of $L$.

\begin{lem} \label{lem2.1} Let $\theta\in\Der_k(B)$ be a 
  $k$-derivation of $B$ such that the following holds:
  \begin{enumerate}
  \item
    the $k$-linear map 
    \[
        \bar{\theta}:\m/\m^2\to k
    \]
    induced by $\theta$ is not zero, 
  \item
    we have $\theta^p=a\theta$, for some $a\in \Frac(B)$. 
  \end{enumerate}
  Then there exists a regular system of parameters $(x_1,\ldots,x_d)$ for $B$
  such that
\[
      \theta(x_1) \equiv 1 \pmod{\m}
\]
and
\[
     \theta(x_2)=\ldots=\theta(x_d)=0.
\]
\end{lem}

\proof Let $\theta$ be as in the statement of the lemma. It follows
immediately from (i) that there exists a system of parameters $x_1,\ldots,x_d$
such that
\begin{equation} \label{eq2.1}
\begin{split}
   \theta(x_1)  &\equiv 1 \pmod{\m},\\
   \theta(x_i)  &\equiv 0 \pmod{\m} \quad\text{for $i>1$.}
\end{split}
\end{equation}
Our strategy is to change the $x_i$ for $i>1$ step by step in order to improve
the last congruence modulo arbitrary powers of $\m$. Since $B$ is complete, we
can take the limit and find parameters $x_i$ such that $\theta(x_i)=0$ for
$i>1$.

Suppose that 
\begin{equation} \label{eq2.2}
   \theta(x_i)  \equiv 0 \pmod{\m^n} \quad\text{for $i>1$,}
\end{equation}
for some $n\geq 1$. We claim that there exist elements
$\Delta_2,\ldots,\Delta_d\in\m^{n+1}$ such that
\begin{equation} \label{eq2.3}
  \theta(\Delta_i)\equiv -\theta(x_i) \pmod{\m^{n+1}}, \qquad i=2,\ldots,d.
\end{equation}
Assuming this claim, we can set
\[
     \tilde{x}_i:=x_i+\Delta_i,\qquad i=2,\ldots,d,
\]
and obtain
\[
    \theta(\tilde{x}_i)\equiv 0 \pmod{\m^{n+1}}, \quad i=2,\ldots,d.
\] 
The lemma then follows by induction and a limit argument. 

To prove the claim we consider the $k$-linear map
\[
     \bar{\theta}_n:\m^{n+1}/\m^{n+2}\to\m^n/\m^{n+1}
\]
induced by $\theta$. We have to show that for $i>1$ the class of $\theta(x_i)$
in $\m^n/\m^{n+1}$ lies in the image of $\bar{\theta}_n$.  A short calculation,
using \eqref{eq2.1} and \eqref{eq2.2}, shows that the image of
$\bar{\theta}_n$ is spanned by the images of the monomials of degree $n$
\[
     x_1^{l_1}x_2^{l_2}\cdots x_d^{l_d}, \qquad l_1+\ldots+l_d=n,
\]
such that
\[
    l_1 \not\equiv -1 \pmod{p}.
\]
In particular, for $n<p-1$ the map $\bar{\theta}_n$ is surjective, and the
claim is true. 

Suppose now that $n\geq p-1$ and fix an index $i>1$. By \eqref{eq2.2} we can
write
\[
  \theta(x_i) \equiv \sum_{l=0}^n \;a_l\,x_1^ly_l \pmod{\m^{n+1}},
\]
where $a_l\in k$ and $y_l$ is a homogenous polynomial of degree $n-l$ in the
variables $x_2,\ldots,x_d$. A straightforward calculation, using \eqref{eq2.1}
and \eqref{eq2.2}, shows that
\begin{equation} \label{eq2.3}
  \theta^p(x_i) \equiv \sum_{l=p-1}^n \;
      l(l-1)\cdots(l-p+2)\,a_l\,x_1^{l-p+1}y_l \pmod{\m^{n-p+2}}.
\end{equation}
On the other hand, it follows from (ii) and \eqref{eq2.2} that
\begin{equation} \label{eq2.4}
     \theta^p(x_i) = a\theta(x_i) \equiv 0 \pmod{\m^n}.
\end{equation}
(Here we have used that the element $a\in\Frac(B)$ occuring in (ii) is of the
form $a=\theta(x_1)^{-1}\theta^p(x_1)$ and therefore lies in $B$.) Combining
\eqref{eq2.3} and \eqref{eq2.4} we find that $a_l=0$ if $l\equiv -1
\pmod{p}$. This shows that the class of $\theta(x_i)$ in $\m^n/\m^{n+1}$ lies
in the image of $\bar{\theta}_n$ and proves the claim. The proof of Lemma
\ref{lem2.1} is now complete.
\Endproof

\begin{cor} \label{cor2.1}
  Let $\theta\in\Der_k(B)$ be as in Lemma \ref{lem2.1}. Then:
  \begin{enumerate}
  \item
    The subring $A:=\Ker(\theta)\subset B$ is regular, and $B/A$ is finite and
    flat of rank $p$. 
  \item
    Suppose that 
    \[
       \theta(f)\in B,
    \]
    for some element $f\in L$. Then there exists $f_0\in B$ such that
    \[
         \theta(f_0)=\theta(f).
    \]
  \item
    Suppose that 
    \[
       \theta(f)/f \in B,
    \]
    for some element $f\in L^\times$. Then there exists $u\in B^\times$ such that
    \[
         \theta(u)/u = \theta(f)/f.
    \]
  \end{enumerate}
\end{cor}

\proof
Let $(x_1,\ldots,x_d)$ be a regular system of parameters as in Lemma
\ref{lem2.1}, and set $h:=\theta(x_1)\in B^\times$. Then $\theta=h\theta_0$,
where 
\[
    \theta_0(x_1)=1,\;\theta_0(x_2)=\ldots=\theta_0(x_d)=0.
\]
Clearly, we have $A=\Ker(\theta)=\Ker(\theta_0)$. It is now easy to see that
\[
    A=k[[x_1^p,x_2,\ldots,x_d]]
\]
is regular and that 
\begin{equation} \label{eq2.7}
    B=A\oplus A\cdot x_1\oplus\ldots\oplus A\cdot x_1^{p-1}
\end{equation}
is finite and flat over $A$ of rank $p$. So (i) is proved. 

For the proof of (ii) we write
\[
      f=a_0+a_1x_1+\ldots+a_{p-1}x_1^{p-1}, 
\]
with $a_i\in K:=\Frac(A)=L\cap A$. Then the hypothesis $\theta(f)\in B$ implies
\[
      \theta_0(f)= a_1+2a_2x_1+\ldots+(p-1)a_{p-1}x_1^{p-2}\in B.
\]
Now \eqref{eq2.7} shows that $a_1,\ldots,a_{p-1}\in A$. Therefore, 
\[
    \theta(f_0)=\theta(f), \qquad\text{with}\;\;
        f_0:=a_1x_1+\ldots+a_{p-1}x_1^{p-1}\in B.
\]
This proves (ii).

Assertion (iii) follows from \cite{Samuel64}, Theorem 2. For the convenience
of the reader, we sketch the argument.

Let $f\in L^\times$ be given such that $\theta(f)/f\in B$.  We claim that
there exists a Weil divisor $\mathfrak{d}$ on $\Spec(A)$ such that
\[
       B\cdot \mathfrak{d} = (f).
\]
(This is obviously a `descent argument' and explains the title of this
section.) 
Assuming the claim we can prove (ii), as follows. By (i), $A$ is regular and in
particular factorial. Therefore, $\mathfrak{d}=(f_0)$, for some $f_0\in
K=\Frac(A)$. So $f=uf_0$ for a unit $u\in B^\times$, and we get
$\theta(f)/f=\theta(u)/u$, as desired. 

To prove the claim, it suffices to show the following: for every prime ideal
$\p\lhd B$ of height one we have
\[
        e_\p \mid v_\p(f).
\]
Here $v_\p$ is the normalized valuation on $L$ associated to $\p$ and $e_\p$
denotes the ramification index of $\p$ in the field extension $L/K$. Since
$e_\p\in\{1,p\}$, we may assume that $n:=v_\p(f)$ is prime to $p$, and we
have to show that $e_\p=1$. 

Let $t\in B_\p$ be a uniformizer for $v_\p$. Then we can write
\[
      f=ut^n, \qquad \text{with}\quad u\in B_\p^\times.
\]
From 
\[
    \theta(f)/f=\theta(u)/u+n\cdot\theta(t)/t\in B
\]
we conclude $\theta(t)/t\in B_\p$. This means that $\theta$ induces a
derivation $\bar{\theta}$ on the residue field $k_\p$. From assumption (i) of
Lemma \ref{lem2.1} it follows that $\bar{\theta}\neq 0$. From
\[
       k_{\p\cap A}\subset k_\p^{\bar{\theta}=0} \subsetneq k_\p
\]
we conclude 
\[
        f_\p:=[k_\p:k_{\p\cap A}] =p.
\]
Now the inequality $p=[L:K]\geq e_\p\cdot f_\p$ implies $e_\p=1$, and the
claim is proved.
\Endproof

\section{Ramification of $p$-cyclic Galois extensions} \label{sec3}

In the situation of Theorem \ref{mainthm}, we obtain a cyclic Galois extension
$L/K$ of degree $p$ by taking fraction fields: $L:=\Frac(B)$,
$K:=\Frac(A)=L^G$.  The present section contains some preliminary
investigation of the ramification of this extension with respect to the
discrete valuation corresponding to the parameter $\pi$. The main point here
is that for $p$-cyclic Galois extensions it is possible to distinguish two
types of wild ramification: {\em total ramification} and {\em fierce
  ramifi\-cation}. Accordingly, the proof of Theorem \ref{mainthm} will be
divided into these two cases.

\subsection{} \label{sec3.1}

Let $L/K$ be a cyclic Galois extension of degree $p$ (where $p$ is, as always,
a prime number). Let $v$ be a discrete valuation on $L$ which is fixed by $G$.  

We choose a generator $\sigma$ of $G$. The valuation rings of $K$ and $L$ with
respect to $v$ are denoted by $\OO_K$ and $\OO_L$, the residue fields by $\Kb$
and $\Lb$. The letter $\pi$ will always denote a uniformizer for $v$ on $L$
(i.e.\ $\pi\OO_L=\p_L$ is the maximal ideal of $\OO_L$). We normalize the
valuation $v$ such that $v(\pi)=1$.

We set
\[
       e_{L/K}:=[v(L^\times):v(K^\times)], \qquad 
       f_{L/K}:=[\Lb:\Kb].
\]
These invariants are related by the fundamental equality
\[
       e_{L/K}\cdot f_{L/K} = [L:K]=p.
\]
See e.g.\ \cite{LangAlgebra}, Corollary XII.6.3. 

\begin{defn} \label{def3.1}
  We classify the ramification behavior of $L/K$ at $v$ as follows: 
  \begin{itemize}
  \item[{\bf I:}] If $e_{L/K}=p$ and $f_{L/K}=1$, we call $L/K$ {\em
    totally ramified}. There are two subcases:
  \begin{itemize}
  \item[(a)]
    If ${\rm char}(\Lb)\neq p$, $L/K$ is called {\em tamely ramified}.
  \item[(b)]
    If ${\rm char}(\Lb)=p$, $L/K$ is called {\em totally wildly ramified}.
  \end{itemize}
  \item[{\bf II:}] 
    Suppose that $e_{L/K}=1$ and $f_{L/K}=p$. There are again two subcases:
    \begin{itemize}
    \item[(a)]
      If $\Lb/\Kb$ is Galois, then $L/K$ is called {\em unramified}.
    \item[(b)]
      If $\Lb/\Kb$ is inseparable, then $L/K$ is called {\em
      fiercely ramified}.
    \end{itemize}
  \end{itemize} 
\end{defn}

The following invariant is useful to distinguish these cases: 
\begin{equation} \label{eq3.9}
     \delta := \min\{\,v(\sigma(x)-x) \mid x\in \OO_L \,\}\geq 0.
\end{equation}

\begin{prop} \label{prop3.1}
\begin{enumerate}
\item
  $L/K$ is unramified if and only if $\delta=0$. 
\item If $L/K$ is tamely ramified, then $\delta=1$. (Note: the converse of
  this implication is actually false: $L/K$ may be fiercely ramified.)
\item
  The extension $L/K$ is totally ramified if and only if
  \begin{equation} \label{eq3.1a}
    \delta=v(\sigma(\pi)-\pi),
  \end{equation}
  for some uniformizer $\pi$. In this case \eqref{eq3.1a} holds for every
  uniformizer $\pi$. 
\end{enumerate}
\end{prop}
 
\proof Assertions (i) and (ii) are classical. See e.g.\ \cite{SerreCL}, IV, \S
1. Assertion (iii) is certainly well known as well. We will nevertheless give a
proof because it yields some useful insight.

We may assume that $\delta>0$. Let $\pi$ be an arbitrary uniformizer. By the
definition of $\delta$, the map
\[
   \tilde{\theta}_\pi:\OO_L\to\Lb, \qquad 
      x\mapsto \,\Big(\,\frac{\sigma(x)-x}{\pi^\delta} \mod{\pi}\,\Big),
\]
is well defined and does not vanish. A short calculation shows that
$\tilde{\theta}_\pi$ is a derivation (i.e.\
$\tilde{\theta}_\pi(x+y)=\tilde{\theta}_\pi(x)+\tilde{\theta}_\pi(y)$ and
$\tilde{\theta}_\pi(xy)=\bar{x}\cdot\tilde{\theta}_\pi(y)+\bar{y}\cdot
\tilde{\theta}_\pi(x)$), and hence we have
\begin{equation} \label{eq3.0}
         \tilde{\theta}_\pi(a\pi) = \bar{a}\cdot\tilde{\theta}_\pi(\pi),
\end{equation}
for all $a\in\OO_L$. 

Suppose that $\tilde{\theta}_\pi(\pi)=0$. This is equivalent to the condition
\begin{equation} \label{eq3.1}
    \delta<v(\sigma(\pi)-\pi).
\end{equation}
Then \eqref{eq3.0} shows that both the equation $\tilde{\theta}_\pi(\pi)=0$
and \eqref{eq3.1} hold in fact for all uniformi\-zers $\pi$. Moreover,
$\tilde{\theta}_\pi$ induces a derivation 
\[
      \theta_\pi:\Lb\to\Lb, \qquad \bar{x}\mapsto \tilde{\theta}_\pi(x)
        \mod{\pi}.
\]
We have $\theta_\pi\neq 0$ by definition of $\delta$. In particular, there
exists a unit $x\in\OO_L^\times$ with $\delta=v(\sigma(x)-x)$. It is also
clear that
\[
       \Kb \subset \Ker(\theta_\pi)\subsetneq\Lb.
\]
Since $[\Lb:\Kb]\leq p$, it follows that $\Kb=\Ker(\theta_\pi)$ and
$[\Lb:\Kb]=p$. In particular, $L/K$ is fiercely ramified (Case II (b)).

On the other hand, if $L/K$ is fiercely ramified, then $e_{L/K}=1$, which
means that there exists a uniformizer $\pi$ in $\OO_K=\OO_L^G$. Then we
obviously have $\tilde{\theta}_\pi(\pi)=0$. All claims made in Proposition
\ref{prop3.1} have now been proved. 
\Endproof

\subsection{} \label{sec3.2}

Suppose that $L/K$ is fiercely ramified (Case II (b)). In this case the proof
of Proposition \ref{prop3.1} yields the following.

\begin{cor} \label{cor3.1}
  Suppose that $L/K$ is fiercely ramified. Then for every uniformizer $\pi$
  the map
  \[
      \theta_\pi:\Lb\to\Lb, \qquad \bar{x}\mapsto 
            \,\Big(\,\frac{\sigma(x)-x}{\pi^\delta} \mod{\pi}\,\Big),
  \]
  is a derivation with the following properties.
  \begin{enumerate}
  \item
    $\theta_\pi\neq 0$,
  \item
    $\Ker(\theta_\pi)=\Kb$,
  \item
    $\theta_\pi^p=a\,\theta_\pi$, for some $a\in\Kb$.
  \end{enumerate}
\end{cor}
  
\proof (i) and (ii) have been proved above. Property (iii) is true for all
derivations of a field $\Lb$ of characteristic $p$ for which
\[
      [\Lb:\Ker(\theta_\pi)]=p.
\]
See e.g.\ \cite{JacobsonGalois}, Chapter IV.8, Exercise 3. 
\Endproof


The following lemma will be used in Section \ref{sec4.4}.

\begin{lem} \label{lem3.2}
  Let $x\in\OO_L$ be such that
  \[
     v(\sigma(x)-x)\geq \delta+n,
  \]
  for some $n\geq 1$. Then there exists 
  $y\in\OO_K$ with
  \[
        y \equiv x \pmod{\p_L^{n}}.
  \]
\end{lem}

\proof We proceed by induction on $n$, starting with $n=1$.  Let $\pi\in\OO_K$
be an invariant uniformizer. The assumption $v(\sigma(x)-x)\geq\delta+1$ means
$\theta_\pi(\bar{x})=0$. It follows from Corollary \ref{cor3.1} (ii) that
$\bar{x}\in\Kb$. So we can take for $y\in\OO_K$ any representative of
$\bar{x}$.

We may hence suppose $n\geq 2$. By induction, there exists an element
$z\in\OO_K$ such that $z\equiv x\pmod{\pi^{n-1}}$. Write
$x=z+a\pi^{n-1}$, with $a\in\OO_L$. Then
\[
     \theta_\pi(\bar{a}) \equiv \frac{\sigma(x)-x}{\pi^{\delta+n-1}}\equiv 0
          \pmod{\pi}.
\]
As before, it follows that $\bar{a}\in\Kb$. Let $b\in\OO_K$ be any lift of
$\bar{a}$ and set $y:=z+b\pi^{n-1}\in\OO_K$. 
\Endproof

\subsection{} \label{sec3.3}

It is well known that the trace map $\Tr_{L/K}:L\to K$ contains interesting
information about the ramification of $L/K$. We will use this only in the case
where $L/K$ is totally wildly ramified.

As before, we fix a uniformizer $\pi\in\OO_L$, $v(\pi)=1$. For any element
$a\in L^\times$ we set
\[
     [a]_\pi:= a\cdot\pi^{-v(a)}\mod{\pi} \in \Lb^\times.
\]

\begin{rem} \label{rem3.1}
  The map $L^\times\to\Lb^\times$, $a\mapsto[a]_\pi$, obeys the following
  rules: 
  \begin{enumerate}
  \item
    $[\pi]_\pi=1$,
  \item
    $[a\cdot b]_\pi = [a]_\pi\cdot[b]_\pi$,
  \item
    $[a+b]_\pi=[a]_\pi+[b]_\pi$, if $v(a)=v(b)$ and $[a]_\pi+[b]_\pi\neq 0$.
  \end{enumerate}
\end{rem}

\begin{lem} \label{lem3.1}
  Suppose that $L/K$ is totally wildly ramified (Case I (b)). Let
  $\pi\in\OO_L$ be a uniformizer (i.e.\ $v(\pi)=1$).
  \begin{enumerate}
  \item
    The trace induces an $\Lb$-linear isomorphism
    \[
        \Trb:\p_L^{-(p-1)(\delta-1)}/\p_L^{-(p-1)(\delta-1)+1} 
             \iso \OO_K/\p_K\cong\Lb.
    \]
  \item
    There exists an element $s_{L/K}^{(\pi)}\in
    \Lb^\times$, uniquely determined by the condition
    \[
        \Tr_{L/K}(x\pi^{-(p-1)(\delta-1)}) \equiv s_{L/K}^{(\pi)}\cdot x \pmod{\pi},
    \]
    for all $x\in\OO_L$.
  \item 
    If $\pi'=u\pi$ is another uniformizer, with $u\in \OO_L^\times$, then
    \[
         s_{L/K}^{(\pi')}= \bar{u}^{-(p-1)(\delta-1)}\cdot s_{L/K}^{(\pi)}.
    \]
  \item
    We have
    \[
        s_{L/K}^{(\pi)} = -[\sigma(\pi)-\pi]_\pi^{p-1}.
    \]
  \end{enumerate}
\end{lem}

\proof (cf.\ \cite{KatoSwan}, Prop.\ 2.4) 
Set $\lambda:=\sigma(\pi)-\pi$. By Proposition \ref{prop3.1} (ii), (iii) and the
assumption on $L/K$ we have $v(\lambda)=\delta\geq 2$ and 
\[
     \sigma(\lambda)\equiv \lambda \pmod{\pi^{\delta+1}}.
\]
Also, for $i=1,\ldots,p-1$ we have
\[
      \sigma^i(\pi)-\pi= \lambda+\sigma(\lambda)+\ldots+\sigma^{i-1}(\lambda).
\]
Using Remark \ref{rem3.1} (iii) and the above, one shows easily that
\begin{equation} \label{eq3.5}
     v(\sigma^i(\pi)-\pi) = v(\lambda)=\delta
\end{equation}
and
\begin{equation} \label{eq3.6}
    [\sigma^i(\pi)-\pi]_\pi = i\cdot[\sigma(\pi)-\pi].
\end{equation}

Let $P\in\OO_L[X]$ denote the minimal polynomial of $\pi$ over $K$. By
\cite{SerreCL}, Lemma III.6.2, we have
\begin{equation} \label{eq3.3}
   \Tr_{L/K}\big(\pi^iP'(\pi)^{-1}\big) = 
      \begin{cases}
        \;\; 0, & i=0,\ldots,p-2, \\
        \;\; 1, & i=p-1.
      \end{cases}
\end{equation}
Since 
\[
      P'(\pi) = \prod_{i=1}^{p-1} \,(\pi-\sigma^i(\pi)),
\]
we get $v(P'(\pi))=(p-1)\delta$ from \eqref{eq3.5}. To prove (i), note that an
element $a\in\p_L^{-(\delta-1)(p-1)}$ can be uniquely written in the form
\[
     a=\sum_{i=0}^{p-1}\,a_i\frac{\pi^i}{P'(\pi)},
\]
with $a_0,\ldots,a_{p-2}\in\p_K$ and $a_{p-1}\in\OO_K$. By \eqref{eq3.3} we
have $\Tr_{L/K}(a)=a_{p-1}$; furthermore, $a_{p-1}\in\p_K$ if and only if
$a\in\p_L^{-(\delta-1)(p-1)+1}$. This proves (i). The Assertions (ii) and
(iii) follow easily. Using \eqref{eq3.3}, \eqref{eq3.6} and the definition of
$s_{L/K}^{(\pi)}$, we get
\[\begin{split}
  s_{L/K}^{(\pi)} & = [P'(\pi)]_\pi 
                  = [\prod_{i=1}^{p-1}\pi-\sigma^i(\pi)]_\pi \\
                & = (p-1)!\cdot [\sigma(\pi)-\pi]^{p-1}_\pi 
                  = -[\sigma(\pi)-\pi]^{p-1}_\pi,
\end{split}\]
proving (iv).
\Endproof

\section{The proof of Theorem \ref{mainthm}} \label{sec4}

In this section we prove our main result, Theorem \ref{mainthm}. The
assumptions made in Section \ref{sec1.1} and \ref{sec1.2} and Conditions (i)
and (ii) of Theorem \ref{mainthm} are in force throughout.

\subsection{} \label{sec4.1}

The ring of invariants $A:=B^G$ is a local, integrally closed domain. By
Assumption (c) of Section \ref{sec1.1}, $A$ is noetherian and
$B$ is a finite $A$-algebra.

Our first step is to show that in the proof of Theorem \ref{mainthm} we may
assume that $B$ is complete. This assumption will be used later in the proof
of Lemma \ref{lem4.4.1} and Lemma \ref{lem4.4.2}.

Let $\Ah$ and $\Bh$ denote the completions of $A$ and $B$, with respect to
their maximal ideals. These are complete local rings with residue field
$k$. The action of $G$ on $B$ extends uniquely to $\Bh$, and we have a
canonical map $\Ah\to\Bh^G$. 

\begin{lem}
  The above map is an isomorphism, $\Ah\cong\Bh^G$.
\end{lem}

\proof By definition of $A$ and Assumption (c) we have an exact sequence of
{\em finite} $A$-modules
\begin{equation} \label{eq4.1.1}
      0 \to A \to B \lpfeil{\sigma-1} B.
\end{equation}
Since $B$ is noetherian we have ${\rm rad}(\m_AB)=\m_B$, so $\Bh$ is the
completion of the $A$-module $B$. Since $A$ is noetherian, the exact sequence
\eqref{eq4.1.1} induces an exact sequence 
\[
      0 \to \Ah \to \Bh \lpfeil{\sigma-1} \Bh,
\]
by \cite{MatsumuraCA}, Thm.\ 54. This proves the lemma. 
\Endproof

By \cite{AtiyahMcdonald}, Prop.\ 11.24, $\Ah$ is regular if and only if $A$ is
regular. So all we have to show is that the action of $G$ on $\Bh$ satisfies
Hypothesis (i) and (ii) of Theorem \ref{mainthm} (if it does for $B$).

Let $x_1,\ldots,x_d$ be a regular system of parameters for $B$; it is also a
regular system of parameters for $\Bh$. Moreover, it is easy to check that
\[
    \I_\sigma=(\sigma(x_1)-x_1,\ldots,\sigma(x_d)-x_d)\lhd B.
\]
Therefore, 
\[
     \hat{\I}_\sigma:=(\sigma(x)-x\mid x\in\Bh)=\Bh\cdot\I_\sigma.
\]
Our claim follows immediately. 

For the rest of this section we may therefore assume that $B$ is a complete
local ring.

\subsection{} 

We set $\Bb:=B/B\pi$; this is a complete local ring with residue field
$k$. Since $\pi$ is a regular parameter of $B$, $\Bb$ is regular. It follows
from Condition (i) of Theorem \ref{mainthm} that the map
\[
    \tilde{\theta}_\pi:B\to\Bb, \qquad 
     x\mapsto\,\Big(\, \frac{\sigma(x)-x}{\pi^\delta} \pmod{\pi}\,\Big)
\]
is a well defined derivation. Moreover, the ideal generated by the image of
$\tilde{\theta}_\pi$ is equal to $\Bb$. 

Let $L:=\Frac(B)$ denote the fraction field of $B$. Let $v:L^\times\to\ZZ$
denote the discrete valuation corresponding to the prime ideal $B\pi\lhd B$,
normalized such that $v(\pi)=1$. So the residue field of $v$ is the fraction
field of $\Bb$, $\Lb=\Frac(\Bb)$. 

Let $K:=L^G$ be the fixed field. Then $L/K$ is a Galois extension with Galois
group $G=\gen{\sigma}\cong\ZZ/p\ZZ$, and $v$ is fixed by $G$. We are therefore
in the situation of Section \ref{sec3}. Clearly, the number $\delta$ from
Theorem \ref{mainthm} is the same as the invariant defined by \eqref{eq3.9}:
\[
    \delta= \min\{\, v(\sigma(x)-x) \mid x\in \OO_L \,\}.
\]
Since $\delta>0$, the extension $L/K$ is ramified, i.e.\
Case II (a) in Definition \ref{def3.1} is excluded. There are three cases
left.

If $\tilde{\theta}_\pi(\pi)\neq 0$, then $L/K$ is totally ramified 
by Proposition \ref{prop3.1} (iii). Moreover, Condition (ii) of Theorem
\ref{mainthm} implies that $\tilde{\theta}_\pi(\pi)\in\Bb^\times$ is a
unit. If $\delta=1$ then $L/K$ is tamely ramified (Case I (a)); otherwise,
$L/K$ is totally wildly ramified (Case I (b)). 

On the other hand, if $\tilde{\theta}_\pi(\pi)=0$, then $L/K$ is fiercely
ramified (Case II (b)). In this case the derivation $\tilde{\theta}_\pi$
induces a nontrivial derivation $\theta_\pi:\Bb\to\Bb$, and the unique
extension of $\theta_\pi$ to $\Lb$ is the derivation from Corollary
\ref{cor3.1}.

After these preliminary remarks, the proof of Theorem \ref{mainthm} is done by a
case-by-case analysis of the three different types of ramification of $v$ in
$L/K$. 

\subsection{Case I (a): $L/K$ is tamely ramified}

Suppose that ${\rm char}(\Lb)\neq p$. By Proposition
\ref{prop3.1} (ii) we have $\delta=1$. Therefore, the morphism
\[
     G \to \Aut(B\pi/B\pi^2)
\]
is injective. It follows that
\[
     \sigma(\pi)\equiv \zeta_p\cdot\pi \pmod{\pi^2},
\]
for some primitive $p$th root of unity $\zeta_p\in\Bb$. But the assumption
$\I_\sigma=(\pi)$ implies that 
\[
     \sigma(\pi)-\pi =u\pi
\]
for a unit $u\in B^\times$. From the above we get 
\[
  u\equiv \zeta_p-1\pmod{\m}. 
\]
In particular, $\zeta_p-1\not\in\m$, which shows that ${\rm char}(k)\neq
p$. Now the hypothesis of Theorem \ref{mainthm} implies that $A$ is regular,
see \cite{Serre67} and Remark \ref{rem1.1}.

\subsection{Case I (b): $L/K$ is totally wildly ramified}

In this case Condition (ii) of Theorem \ref{mainthm} implies
\begin{equation} \label{eq4.3}
        \sigma(\pi)-\pi = u\pi^\delta,
\end{equation}
where $u\in B^\times$ is a unit and $\delta\geq 2$. Let $\p:=A\cap B\pi\lhd
A$. This is a prime ideal of height one.

\begin{lem}
  The element
  \[
       \lambda:=N_{L/K}(\pi) =\prod_{i=0}^{p-1} \sigma^i(\pi) \in A
  \]
  generates $\p$. In particular, $\p$ is a principal ideal.
\end{lem}

\proof It is obvious that $\lambda\in\p$. Let $\alpha\in\p$ be an arbitrary
element. To prove that $\p=(\lambda)$, it suffices to show that $\lambda\mid
\alpha$ in $B$.  

We have $\pi\mid\alpha$ by definition of $\p$. This means that
$v(\alpha)>0$. Since $\alpha\in A\subset K$ and $e_{L/K}=p$ by assumption, we
actually have $v(\alpha)\geq p$. This implies
$\pi^p\mid \alpha$. 

By \eqref{eq4.3} we have
\[
     \sigma(\pi) = \pi\,(1+u\pi^{\delta-1}) \sim \pi.
\]
Here $a\sim b$ means: $b=va$ for a unit $v\in B^\times$.  It follows that
$\sigma^i(\pi)\sim\pi$ for all $i$ and hence $\lambda\sim\pi^p$. Therefore,
$\lambda|\alpha$, as desired, and the lemma is proved.
\Endproof

Let $\Ab:=A/\p$. We consider $\Ab$ as a subring of $\Bb$ via the natural
embedding $\Ab\inj\Bb=B/B\pi$.

\begin{lem}
  We actually have $\Ab=\Bb$. In particular, $\Ab$ is regular. 
\end{lem}

\proof
Clearly, the rings $\Ab$ and $\Bb$ have the same fraction field
$\Lb$. Consider the $\Lb$-linear map 
\[
   \Trb:\Lb\to\Lb, \qquad \bar{x}\mapsto 
      \,\big(\, \Tr_{L/K}(x\pi^{-(p-1)(\delta-1)}) \mod{\pi}\big),
\]
studied in Section \ref{sec3.3}. We claim that $\Trb(\Bb)\subset\Ab$. Indeed,
if $x$ is an element of $B$, then the proof of Lemma \ref{lem3.1} shows that
\[
     v(\Tr_{L/K}(x\pi^{-(p-1)(\delta-1)}))\geq 0. 
\]
Therefore, $\Tr_{L/K}(x\pi^{-(p-1)(\delta-1)})\in (B[\pi^{-1}]\cap\OO_L)^G=B^G=A$
by factoriality of $B$, proving the claim.

It follows from \eqref{eq4.3} and Lemma \ref{lem3.1} (ii), (iv) that
\[
     \Trb(\xb) = -\bar{u}^{p-1}\xb,
\]
where $\bar{u}\in\Bb^\times$ is the image of $u$ in $\Bb$. We conclude
$-\bar{u}^{p-1}\Bb=\Bb\subset\Ab$. So $\Ab=\Bb$, and the lemma is proved.
\Endproof

Using the two lemmas it is now easy to see that $A$ is regular: if
$\xb_1,\ldots,\xb_{d-1}$ is a regular sequence of parameters for $\Ab$ and
$x_i\in A$ is a lift of $\xb_i$, then $(\lambda,x_1,\ldots,x_{d-1})$ is a
regular system of parameters for $A$.

\subsection{Case II (b): $L/K$ is fiercely ramified} \label{sec4.4}

Let $\theta=\theta_\pi\in\Der(\Lb)$ be the derivation from  Corollary
\ref{cor3.1}. It follows from the definition of $\theta$ that 
\[
      \theta(\Bb)\subset\Bb.
\]
Also, Condition (i) of Theorem \ref{mainthm} implies that there exists an
element $x\in B$ such that
\[
      \frac{\sigma(x)-x}{\pi^\delta}\in B^\times.
\]
This shows that the $k$-linear map 
\[
      \bar{\theta}:\m_{\Bb}/\m_{\Bb}^2\to k.
\]
induced by $\theta$ is not zero, $\bar{\theta}\neq 0$. Since
$\theta^p=a\,\theta$ for some element $a\in\Kb$ (Corollary \ref{cor3.1}
(iii)), the derivation $\theta$ satisfies the hypotheses of Lemma \ref{lem2.1}
(with respect to the ring $\Bb$). Therefore, we may (and will) use Corollary
\ref{cor2.1} in the proof of the following two lemmas. 

\begin{lem} \label{lem4.4.1}
  There exists an element $\lambda\in A$ of the form
  $\lambda=u\pi$, with $u\in B^\times$. 
\end{lem}

\proof By induction, we will construct a sequence of elements
$\lambda_0,\lambda_1,\ldots\in B$ of the form $\lambda_n=u_n\pi$, $u_n\in
B^\times$, such that:
\begin{itemize}
\item[(a)]
  $\lambda_n\equiv \lambda_{n-1}\pmod{\pi^n}$,
\item[(b)]
  $v(\sigma(\lambda_n)-\lambda_n)\geq\delta+n+1$,
\end{itemize}
for all $n\geq 0$. The limit $\lambda:=\lim_n\lambda_n$ is then an element
$\lambda\in A=B^G$ of the form $\lambda=u\pi$, with $u\in B^\times$. 

We start the induction at $n=0$ by setting $\lambda_0:=\pi$. Then
Condition (a) is empty and (b) is true by Proposition \ref{prop3.1} (iii).

The case $n=1$ is special. The uniformizer $\lambda_1$ we are looking for is
of the form $\lambda_1=u_1\pi$, with $u_1\in B^\times$. Condition (a) is still
empty, and Condition (b) is equivalent to 
\[
    0\equiv \frac{\sigma(\lambda_1)-\lambda_1}{\pi^{\delta+1}}
     \equiv u_1w +\theta_\pi(\bar{u}_1) \pmod{\pi},
\]
where $w:=(\sigma(\pi)-\pi)\pi^{-\delta-1}\in B$. So we have to find an
element $\bar{u}_1\in\Bb^\times$ with
\begin{equation} \label{eq4.11}
  \frac{\theta_\pi(\bar{u}_1)}{\bar{u}_1} = -\bar{w}.
\end{equation}
By Lemma \ref{lem3.2}, there exists an element $\mu\in\OO_K$
with $\mu\equiv\pi\pmod{\p_L^2}$. If we write $\mu=v\pi$, with
$v\in\OO_L^\times$, then the same calculation that resulted in \eqref{eq4.11}
shows that
\[
   \frac{\theta_\pi(\bar{v})}{\bar{v}} = -\bar{w}\in\Bb.
\]
Now it follows from Corollary \ref{cor2.1} (iii) that there exists
$\bar{u}_1\in\Bb^\times$ such that \eqref{eq4.11} holds. Now the case $n=1$ of
the lemma is proved. 

We may hence assume that $n\geq 2$ and that $\lambda_0,\ldots,\lambda_{n-1}$
have already been constructed. Since $\lambda_{n-1}\sim\pi$, we may assume
that $\lambda_{n-1}=\pi$. In particular,
\begin{equation} \label{eq4.12}
    v(\sigma(\pi)-\pi) \geq \delta+n.
\end{equation}
The element $\lambda_n$ we are looking for is of the form
\[
        \lambda_n=\pi(1+a\pi^n), 
\]
for some $a\in B$. Using \eqref{eq4.12} we get
\[
   \frac{\sigma(\lambda_n)-\lambda_n}{\pi^{\delta+n}}
     \equiv c + \theta_{\pi}(\bar{a}) \pmod{\pi},
\]
where $c:=(\sigma(\pi)-\pi)\pi^{-\delta-n}\in B$. This means that we have to
find an element $\bar{a}\in\Bb$ such that
\begin{equation} \label{eq4.13}
     \theta_\pi(\bar{a}) = -\bar{c}.
\end{equation}
By \eqref{eq4.12} and Lemma \ref{lem3.2}, there exists an element
$\mu\in\OO_K$ with $\mu\equiv\pi\pmod{\p_L^n}$. Write $\mu=\pi(1+b\pi^n)$,
with $b\in\OO_L$. Then the same calculation from which we deduced
\eqref{eq4.13} yields
\[
    \theta_\pi(\bar{b})=-\bar{c}.
\]
Using Corollary \ref{cor2.1} (ii) we find an element $\bar{a}\in\Bb$
satisfying \eqref{eq4.13}. Now the proof of the lemma is complete.
\Endproof

Let $\lambda\in A$ be as in the lemma.  We get
a natural embedding 
\[
    \Ab:=A/A\lambda \inj \Bb=B/B\pi.
\]
It is clear that $\Ab$ is a complete local ring with residue field $k$ and
quotient field $\Kb$. 

\begin{lem} \label{lem4.4.2}
  We have $\Ab=\Bb\cap\Kb$. 
\end{lem}

\proof
The proof is very similar to the proof of the previous lemma. The inclusion
$\bar{A}\subset\Bb\cap\Kb$ is clear. To prove the converse, suppose we have
an element $x\in B$ such that $\bar{x}\in\Kb$. We will inductively construct a
sequence $y_0,y_1,\ldots\in B$ such that
\begin{itemize}
\item[(a)]
  $x\equiv y_n \pmod{\pi}$,
\item[(b)]
  $y_n\equiv y_{n-1} \pmod{\pi^n}$,
\item[(c)]
  $v(\sigma(y_n)-y_n)\geq \delta+n+1$,
\end{itemize}
for all $n\geq 0$. The limit $y:=\lim_ny_n$ is then an element $y\in A=B^G$
with $x\equiv y\pmod{\pi}$. 

For $n=0$ we set $y_0:=x$. Then (a) is clear, (b) is empty and (c) follows
from $\bar{x}\in\Kb$. So we may assume $n\geq 1$. The element $y_n$ we are
looking for is of the form $y_n=y_{n-1}+a\pi^n$, with $a\in B$. The crucial
Condition (c) is equivalent to 
\[
   \frac{\sigma(y_n)-y_n}{\pi^{\delta+n}} \equiv 
                           c +\theta_\pi(\bar{a}) \equiv 0 \pmod{\pi},
\]
where $c:=(\sigma(y_{n-1})-y_{n-1})\pi^{-\delta-n}\in B$. So we have to find
$\bar{a}\in\Bb$ such that
\begin{equation} \label{eq4.14}
       \theta_\pi(\bar{a}) = -\bar{c}. 
\end{equation}
By Lemma \ref{lem3.2}, there does exist an element $z\in\OO_L$ satisfying the
properties required for $y_n$. If we write $z=y_{n-1}+b\pi^n$, with
$b\in\OO_L$, then we get $\theta_\pi(\bar{b})=-\bar{c}\in\Bb$. Now Corollary
\ref{cor2.1} (ii) shows that $\bar{a}\in\Bb$ satisfying \eqref{eq4.14}
exists, and we are done.
\Endproof

The lemma shows that $\Ab=\Ker(\theta_\pi|_{\Bb})$. By Corollary \ref{cor2.1}
(i) it follows that $\Ab$ is regular. With the same argument used at the end
of the previous subsection we conclude that $A$ is regular. Now Theorem
\ref{mainthm} is proved.
\Endproof

\end{document}